\input epsf

\magnification=1200


\hsize=125mm             
\vsize=195mm             
\parskip=0pt plus 1pt    
\clubpenalty=10000       
\widowpenalty=10000      
\frenchspacing           
\parindent=8mm           

\let\txtf=\textfont
\let\scrf=\scriptfont
\let\sscf=\scriptscriptfont
\font\frtnrm =cmr12 at 14pt
\font\tenrm  =cmr10
\font\ninerm =cmr9
\font\sevenrm=cmr7
\font\fiverm =cmr5

\txtf0=\tenrm
\scrf0=\sevenrm
\sscf0=\fiverm

\def\rm{\fam0 \tenrm}

\font\frtnmi =cmmi12 at 14pt
\font\tenmi  =cmmi10
\font\ninemi =cmmi9
\font\sevenmi=cmmi7
\font\fivemi =cmmi5

\txtf1=\tenmi
\scrf1=\sevenmi
\sscf1=\fivemi

 \def\oldstyle{\fam1 \tenmi}

\font\tensy  =cmsy10
\font\ninesy =cmsy9
\font\sevensy=cmsy7
\font\fivesy =cmsy5

\txtf2=\tensy
\scrf2=\sevensy
\sscf2=\fivesy

\font\tenit  =cmti10
\font\nineit =cmti9
\font\sevenit=cmti7
\font\fiveit =cmti7 at 5pt

\txtf\itfam=\tenit
\scrf\itfam=\sevenit
\sscf\itfam=\fiveit

\def\it{\fam\itfam\tenit}
\font\tenbf  =cmb10
\font\ninebf =cmb10 at  9pt
\font\sevenbf=cmb10 at  7pt
\font\fivebf =cmb10 at  5pt

\txtf\bffam=\tenbf
\scrf\bffam=\sevenbf
\sscf\bffam=\fivebf

\def\bf{\fam\bffam\tenbf}

\newfam\msbfam       
\font\tenmsb  =msbm10
\font\sevenmsb=msbm7
\font\fivemsb =msbm5

\txtf\msbfam=\tenmsb
\scrf\msbfam=\sevenmsb
\sscf\msbfam=\fivemsb

\def\msb{\fam\msbfam\tenmsb}
\def\Bbb#1{{\msb #1}}

\def\RR{{\Bbb R}}

\newfam\scfam
\font\tensc  =cmcsc10

\txtf\scfam=\tensc

\def\sc{\fam\scfam\tensc}


\def\frtnmath{%
\txtf0=\frtnrm        
\txtf1=\frtnmi         
}

\def\frtnpoint{%
\baselineskip=16.8pt plus.5pt minus.5pt%
\def\rm{\fam0 \frtnrm}%
\def\oldstyle{\fam1 \frtnmi}%
\everymath{\frtnmath}%
\everyhbox{\frtnrm}%
\frtnrm }


\def\ninemath{%
\txtf0=\ninerm        
\txtf1=\ninemi        
\txtf2=\ninesy        
\txtf\itfam=\nineit      
\txtf\bffam=\ninebf      
}

\def\ninepoint{%
\baselineskip=10.8pt plus.1pt minus.1pt%
\def\rm{\fam0 \ninerm}%
\def\oldstyle{\fam1 \ninemi}%
\def\it{\fam\itfam\nineit}%
\def\bf{\fam\bffam\ninebf}%
\everymath{\ninemath}%
\everyhbox{\ninerm}%
\ninerm }



\def\text#1{\hbox{\rm #1}}

\def\cite#1{{\uppercase{#1}}}
\def\ref#1{{\uppercase{#1}}}
\def\label#1{{\uppercase{#1}}}
\def\br{\hfill\break} 



\def\topmatter{\null\firstpagetrue\vskip\bigskipamount}
\def\endtopmatter{\vskip2\bigskipamount}

\def\title#1{%
\vbox{\raggedright\frtnpoint
\noindent #1\par}
\vskip 2\bigskipamount}        


\def\shorttitle#1{\rightheadtext={#1}}              

\newif\ifThanks
\global\Thanksfalse

\def\author#1{\begingroup\raggedright
\noindent{\sc #1\ifThanks$^*$\else\fi}\endgroup
\leftheadtext={#1}\vskip \bigskipamount}

\def\endabstract{\endgroup}

\long\def\abstract#1\endabstract{\par
\begingroup\ninepoint\narrower
\noindent{\sc Abstract.\enspace}#1%
\vskip\bigskipamount\endabstract}

\def\section#1#2{\bigbreak\bigskip\begingroup\raggedright
\noindent{\bf #1.\quad #2}\nobreak
\medskip\endgroup\noindent\ignorespaces}

\def\proclaim#1{\medbreak\noindent{\sc #1.\enspace}\begingroup
\it\ignorespaces}
\def\endproclaim{\endgroup\bigbreak}

\def\remark#1{\medbreak\noindent{\sc Remark \enspace}
\begingroup\ignorespaces}
\def\endremark{\endgroup\bigbreak}


\def\qed{$\mathord{\vbox{\hrule\hbox{\vrule
\hskip5pt\vrule height5pt\vrule}\hrule}}$}

\def\demo#1{\medbreak\noindent{{#1}.\enspace}\ignorespaces}
\def\enddemo{\penalty-100\null\hfill\qed\bigbreak}

\newdimen\EZ

\EZ=.5\parindent

\newbox\itembox

\newdimen\ITEM
\newdimen\ITEMORG
\newdimen\ITEMX
\newdimen\BUEXE

\def\iteml#1#2#3{\par\ITEM=#2\EZ\ITEMX=#1\EZ\BUEXE=\ITEM
\advance\BUEXE by-\ITEMX\hangindent\ITEM
\noindent\leavevmode\hskip\ITEM\llap{\hbox
to\BUEXE{#3\hfil$\,$}}%
\ignorespaces}


\newif\iffirstpage\newtoks\righthead
\newtoks\lefthead
\newtoks\rightheadtext
\newtoks\leftheadtext
\righthead={\ninepoint\rm\hfill{\the\rightheadtext}\hfill\llap{\folio}}
\lefthead={\ninepoint\rm\rlap{\folio}\hfill{\the\leftheadtext}\hfill}
\headline={\iffirstpage\hfill\else
\ifodd\pageno\the\righthead\else\the\lefthead\fi\fi}
\footline={\iffirstpage\hfill\global\firstpagefalse\else\hfill\fi}

\leftheadtext={}
\rightheadtext={}


\def\Refs{\bigbreak\bigskip\noindent{\bf References}\medskip
\begingroup\ninepoint\parindent=40pt}
\def\endRefs{\par\endgroup}
\def\endref{}

\def\ref{\par}
\def\key#1{\item{\hbox to 30pt{[#1]\hfill}}}
\def\by{}


\def\cline#1{\leftline{\hfill#1\hfill}}

\def\bR{\Bbb R}
\def\bC{\Bbb C}

\def\log{\mathop{\hbox{\rm log}}\nolimits}

\input epsf

\topmatter
\title{Planar trees, slalom curves and hyperbolic knots}
\author{Norbert A'Campo}
\shorttitle{Trees and knots.}
\endtopmatter

\noindent \S 1. {\bf Introduction}\par
\noindent An embedded tree $B$ in the unit disk $D$, such that the 
intersection $B \cap \partial{D}$ consists of one  
terminal vertex $r$ of $B,$ is called a 
rooted planar tree. For a rooted planar tree $B$ there exists an immersed copy
$P_B \subset D$ of the
interval $[0,1]$ with the following properties:\br
(i) The immersion is relative, i.e. the endpoints are embedded
in $\partial{D}$.\br
(ii) The immersion is generic, i.e. there are
only transversal crossing points, only the endpoints lie
on $\partial{D}$ and the immersion is transversal to $\partial{D}$.\br
(iii) The double points of $P_B$ lie in the interior of the edges of
$B,$ such that the local branches are transversal to the edge of $B$.\br
(iv) Each connected component of $D \setminus P_B$ contains
exactly one vertex of
$B$.\br
(v) The only intersection points of $P_B$
with $B$ are the double points of $P_B.$\br

\noindent
The immersed curve $P_B$ is well defined up to
regular relative isotopy
and is called the slalom curve or slalom divide  of the rooted
planar tree $B$, see Fig. $1,2,3$.

\midinsert
\cline{\epsffile{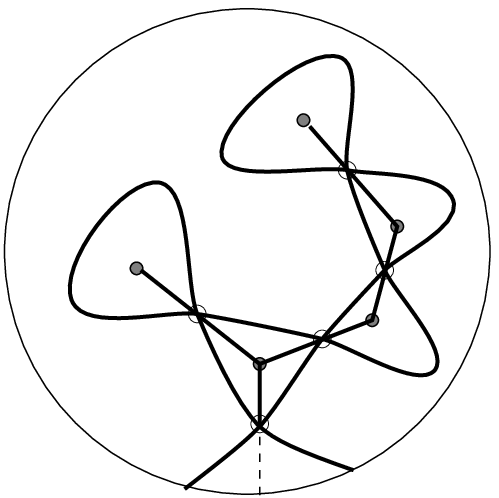}}
\centerline { Figure $1$:
Rooted planar tree,  its Dynkin diagram $E_{10}$ and
slalom.}
\medskip
\endinsert

\noindent 
The slalom curve $P_B$ is a divide to which
corresponds a classical knot $K_B$ in $S^3$, which we call a slalom knot. The
complement of the slalom knot $K_B$ admits
a fibration over the circle $S^1,$ see [AC4] and Section $2$ for basic
definitions and properties. The Dynkin diagram
$\Delta_B$ of the
divide $P_B$ is deduced from the rooted tree $B$ as follows:
First make a new tree $B'$ by subdividing each edge of $B$ with a new
vertex, which is placed at  the crossing point of $P_B$ on the edge; next,
remove from $B'$ the root vertex $r$ and the terminal edge of
$B'$
pointing to $r.$ In Fig. $1$ the tree $B$ has the shape of the classical
Dynkin diagram $D_6$ but the Dynkin diagram  $\Delta_B$ of $P_B$ has $10$
vertices which we can denote by $E_{10}$.
The Dynkin diagram $\Delta_B$ of a rooted tree $B$ is a bicolored rooted
tree with an embedding in the plane.
The root is the new vertex which lies on the edge of
$B$ originating from the root
point of $B$ and the bicoloring is such that the new
vertices are of the same color. Moreover, the Dynkin
diagram $\Delta_B$ has the
property that the terminal vertices of $\Delta_B$ different
from the root, are never
new. The
purpose of this paper is to prove the following theorem.

\proclaim{Theorem 1}
Let $B$ be a rooted tree. The complement of the slalom knot $K_B$
admits a complete hyperbolic metric of finite
volume,
if and only if the Dynkin diagram $\Delta_B$ is 
neither the diagram $A_{2k},\ 1 \leq k,$ nor
the
diagram $E_6$ or $E_8$.
\endproclaim

If the Dynkin diagram $\Delta_B$ is among  $A_{2k},\ 1 \leq k,\ E_6,\ E_8$, the knot
$K_B$ is the torus knot $(2,2k+1),(3,4)$ or $(3,5)$ and  appears
as local knot of a simple
plane curve singularity [AC1]; the monodromy diffeomorphism (with free boundary)
of the knot $K_B$ can be chosen to be
of finite order in those cases and its complement does not carry a
complete hyperbolic metric. We only need to prove the if part of the theorem.

From the above theorem we get many examples of hyperbolic fibered knots,
whose monodromy diffeo\-morphism  and gor\-dian number are known explicitly. 
The monodromy diffeomorphism of a slalom knot  
can be realized as the
product of right Dehn twists of a system of simple closed curves on 
the fiber surface, such that the union of the curves is a spline in the fiber
surface and the dual graph of the system is the Dynkin diagram of the 
rooted tree;
the gordian number of a slalom knot equals the number 
of crossings of the slalom divide [AC4]. 
We call (see section 3) the isotopy class of the 
monodromy diffeomorphism of the slalom knot of a rooted tree the 
Coxeter diffeomorphism of the Dynkin diagram of the rooted tree. 
It follows from Theorem $1$ that a Coxeter diffeomorphism of the Dynkin
diagram of a rooted tree 
is pseudo-anosov, if and only if the Dynkin diagram  
is not a classical Dynkin diagram
(see Theorem $3$). We do not know the lattice in ${\text iso}(H^3)=PSL(2,\bC)$
of the hyperbolic uniformization for the complement of the hyperbolic slalom
knots $K_B$. I like to thank Makoto Sakuma for explaining to me his joint 
work
with Jeff Weeks on hyperberbolic 2-bridge links [S-W], which indicates a road
leading to a description of  the uniformization lattice and the canonical 
decomposition in
ideal
hyperbolic simplices of the complement of hyperbolic 
slalom knots. 

\noindent \S 2. {\bf Divide and knot  of a planar rooted tree.}\par
\noindent Let $B$ be a rooted planar tree in the unit disk 
$D \subset \bR^2$ and
let $P_B$ be its
divide. The knot $K_B$ of the tree $B$ is the 
knot of its divide $P_B$ (see
[AC3-4]), i.e.
$$
K_B:=\{(x,u) \in T(P_B) \ \mid \ \|(x,u)\|=1\} \subset S(T(\bR^2))=S^3
$$
where $T(P_B) \subset T(\bR^2)$ is the subspace of tangent vectors to the 
divide  $P_B$ in the space of tangent vectors to the plane $\bR^2.$
A tangent vector of the plane $(x,u) \in T(\bR^2)=\bR^2 \times \bR^2,$
is represented by its foot
$x \in \bR^2$  and its  linear part $u \in T_x(\bR^2)=\bR^2$. 
The norm $\|(x,u)\|$ is the usual 
euclidean norm of $\bR^4.$ In Fig. $4$ is shown a computer drawing of the 
knot of the
divide Lys (see Fig. $3$). This knot can 
be presented with $11$ crossings and its
gordian number equals the number of crossing points of the divide, i.e. 
$4$. 
Since a slalom divide is connected, the complement of the 
knot $K_B$ of a rooted tree
fibers over the circle [AC4]. A 
model for the fiber surface  and monodromy diffeomorphism 
can be read
by a graphical algorithm from the divide $P_B$ as follows:  replace each
crossing point of $P_B$ by a square, which has its vertices on the local
branches of $P_B$ at the crossing point, and get a trivalent graph $\Gamma$
embedded in the disk $D;$ the
fiber is diffeomorphic to the interior of the surface with boundary 
$F$ obtained from a  thickening of the  graph 
$\Gamma.$ The thickening corresponds to the
cyclic ordering of the edges of $\Gamma$ at each vertex of $\Gamma$, which 
alternatingly agrees or disagrees with an orientation of the ambient plane. 
The graph
$\Gamma$ has only circuits of even length, so the alternating cyclic
ordering of the edges at the vertices of $\Gamma$ exists. 
For each of its squares and for each region of 
the divide $P_B$ the graph $\Gamma$ has  a circuit, which 
surrounds the square or region. To these 
circuits of $\Gamma$
correspond simple closed curves on the surface $F.$ The monodromy $T$ is the
product
of the right Dehn twists along those closed curves. 

\midinsert
\cline{\epsffile{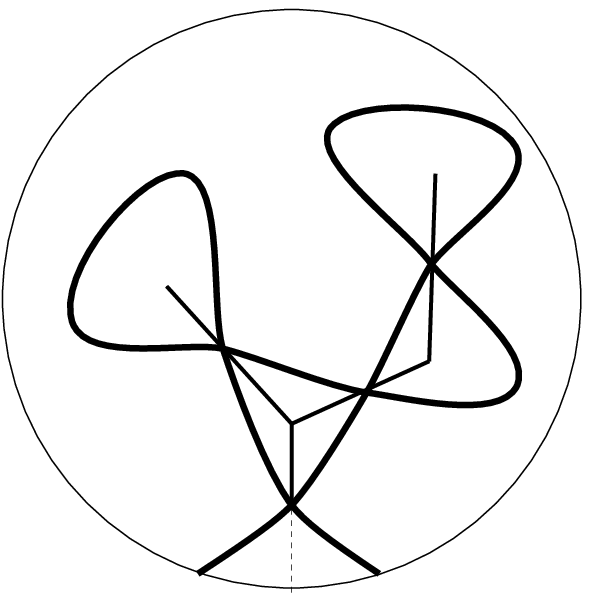}}
\medskip
\centerline{Figure $2$: The slalom $E_{8}$.}
\endinsert

This product is well
defined up to conjugacy in the relative mapping class group of the surface 
$(F,\partial F)$ since
the non-commutation graph of this set of Dehn twists is precisely the Dynkin
diagram $\Delta_B$, which is a tree 
([B], Fascicule XXXIV, Chap. $4$, par. $6$, lemme $1$  ). The
graph $\Gamma$ with its cyclic orientation of 
the edges at the vertices allows us
to give a
combinatorial description of the diffeomorphism $T$ of the surface $F$, which
can be used as input to the Bestvina-Handel algorithm, [B-H] see 
also [L] . In practice, we use
the Cayley code for rooted trees, see [S-Wh], and a 
maple program to deduce from the Cayley
code the combinatorial description of $T$, which was finally the input to the
program TRAINS, written by Tobi Hall, doing the Bestvina-Handel algorithm. 
This way we
get extra stimulating evidence for Theorem $1$ and $3$. 
I would like to thank Tobi Hall for allowing me to use his program TRAINS.

\goodbreak
\noindent \S 3. {\bf Conway spheres and Bonahon-Siebenmann decomposition 
for slalom knots}

Let $B$ be a rooted tree with slalom divide $P_B \subset D$ and 
slalom knot $K_B$. 
Let $f_B:D \to \RR$ be a morse function for the divide $P_B$ as in the 
proof of the fibration theorem of [AC4], i.e. a generic
$C^{\infty}$ function,
such that $P_B$ is its $0$-level and that each interior region 
has exactly one
non-degenerate minimum and that each region which
meets the boundary
has  exactly one non-degenerate maximum or minimum on the
intersection of the region with $\partial D.$ 
The underlying tree of the
slalom divide can be
reconstructed up to isotopy as the closure of the union of 
the gradient lines of $f_B$,
which lie in $\{f_B < 0 \}$ and which contain a saddle point in their closure.
The singular gradient lines $L$ of $f_B$ in $\{f_B > 0 \}$ give enough 
Conway spheres to build the Bonahon-Siebenmann decomposition [B-S] of the knot
$K_B$, see [K]. For a singular gradient line $L$ of $f_B$ in $\{f_B > 0 \}$ 
we define 
$C(L):=\{(x,u)\in T(D) \mid x \in L,\, ||x||^2+ ||u||^2=1\}$. 
Observe that such a gradient
line passes through a saddle point of $f_B$ and both end points of $L$ are on 
$\partial{D}$. It follows that $C(L)$ is a smooth embedded $2$-sphere 
in $S^3$.
Each sphere $C(L)$ is invariant under the involution $(x,u) \mapsto (x,-u)$.

We state without proof:

\proclaim{Theorem 2}
Let $B$ be a rooted tree and $f_B$ its morse function. The spheres $C(L)$ of
the singular gradient $L$ lines of $f_B$ in $\{f_B > 0 \}$ are Conway 
spheres for
the slalom knot $K_B$. 
The spheres $C(L)$ which correspond to edges of the tree $B$,
with at least one endpoint of valency $\geq 3$ or equal to the root vertex,
give the Bonahon-Siebenmann decomposition of the slalom knot.
\endproclaim

\midinsert
\cline{\epsffile{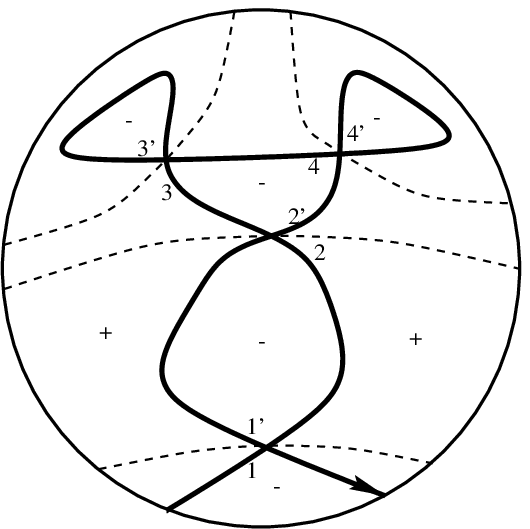}} 
\medskip
\centerline{Figure $2$: The slalom divide of the tree $[0,1,2,2]$ and
singular gradient lines.}
\endinsert

We like to mention here that the 
knot of the slalom divide of the tree $[0,1,2,2]$  is the knot
$10_{139}$ of the table of Rolfson's book [R], which is equivalent  to 
the Montesinos knot
$M(1,(3,1),(3,1),(4,1))$, see [Ka]. 
The gordian number of $10_{139}$ is shown to be $4$ by Tomomi Kawamura [Kaw]. 

I am grateful to Mikami Hirasawa for explaining to me his  
method of  constructing a knot diagram
for slalom knots directly from the slalom divide.
He first doubles the slalom divide and 
then changes according to local rules
the doubled divide to a knot diagram. It follows from
his construction that slalom knots are arboricents knots in the sense of
Bonahon and Siebenmann [B-S]. The notation as arboricent knot for the 
slalom knot $P_B$ is the Dynkin diagram $\Delta_B$  with the weighting $2$
at each vertex. 

\noindent \S 3. {\bf Trees, forms, plumbings.}\par

\noindent Let $A$ be a tree with vertex set $\{v_1, v_2, \dots , v_n\}.$
We will choose the
numbering of the vertices such that for some $m,\ 1 \leq m \leq n,\ $ there
are no
pairs of vertices  $v_i$ and $v_j$ connected by an edge of $A$ with
$i\leq m$ and $j \leq m$ or with $m <i$ and $m < j.$ The chosen numbering
corresponds to a bicoloring of the vertices of the tree.
The real vectorspace
$V_A$
generated by the set of vertices of $A$ carries a quadratic form $q_A,$ whose
matrix is 
$q_A(v_i,v_i)=-2,\ 1\leq i \leq n,\ $ and $q_A(v_i,v_j)=1$, if and only if, the
vertices $v_i$ and $v_j$ are connected by an edge of $A.$ 
To each vertex $v_i$ corresponds an isometry $R_i$ of $(V_A,q_A)$

$$
R_i(v_j):=v_j+q_A(v_i,v_j)v_i,
$$

\noindent
which is a reflection.
Since the non-commutation graph of the set $\{R_1,R_2, \dots ,R_n\}$ 
is a tree 
the product of the reflections $R_i$ does up to conjugacy not depend 
on the
order in which the product is evaluated [B] and 
is called the Coxeter element $C_A$ 
of the tree $A.$ The vector space $V_A$ also 
carries a skew form $sq_A,$ whose
matrix is $sq_A(v_i,v_j)=1$ or $sq_A(v_i,v_j)=-1$ 
if and only if the vertices $v_i$
and $v_j$ are connected by an edge of $A$. 
If $i \leq m$ then $sq_A(v_i,v_j)=1$ else if $j \leq m$ then 
$sq_A(v_i,v_j)=-1$. 
To each vertex $v_i$ corresponds an
endomorphism $T_i$ of $(V_A,sq_A)$

$$
T_i(v_j):=v_j+sq_A(v_i,v_j)v_j,
$$

\noindent
which is a transvection. The product of the transvections $T_i,$ 
equals to $-C_A,$ and we call its conjugacy class in the group of the form
$sq_A$, well defined by [B],  
the skew Coxeter element $sC_A$ of the tree $A.$

From [AC2] we recall the following (see also [H]). 
If the tree $A$ is  
not among the diagrams 
$A_{k},\ D_{k+3},\ \tilde D_{k+3}, 
\ 1 \leq k,\ E_6,\ \tilde E_6,\ E_7,\ \tilde E_7,\ E_8,\ \tilde E_8$, 
the endomorphism
$C_A$ has a real eigenvalue $\lambda_{max}>1$ with 
multiplicity $1$, such that for any eigenvalue
$\lambda$ of $C_A$ we have $|\lambda| < |\lambda_{max}|$ unless
$\lambda=\lambda_{max}$. We call $\lambda_{max}$ the dominating eigenvalue of
$C_A$ and $-\lambda_{max}$ the dominating eigenvalue of
$sC_A$.

Let $A$ be a planar tree with vertex
set  $\{v_1, v_2, \dots ,v_n\}$. To
the
tree $A$ corresponds a surface $S_A$ by the following plumbing.
First realize the planar tree $A$ by a planar circle packing 
with small overlappings. 
Each vertex $v_i$ is
represented by an oriented 
circle $c_i$. As orientation we choose the counterclockwise orientation. 
The circles $c_i,c_j$ are disjoint if the vertices $v_i$ and
$v_j$ are not connected in $B$ and touch each other from the outside with a
small overlap, if $v_i$ and $v_j$ are connected in $A.$   Let $C_i$ be  a
tubular
neighborhood in the plane of $c_i,$ which is an oriented cylinder. 
The surface $S_A$ is obtained by plumbing the cylinders $C_i$ and $C_j$ at one
of the intersection points of $c_i$ and $c_j$ and making an overcrossing 
at the other intersection point,
if the vertices $v_i$ and $v_j$ are
connected in $A.$ The choice at which intersection point the plumbing takes
place, is made such that on the surface $S_A$ the cycles $c_i$ and $c_j$ have
the intersection number $sq_A(v_i,v_j).$ The surface $S_A$ is naturally
immersed
in the plane. Let $D_i$ be the right Dehn twist with core the curve $c_i$ of 
$S_A.$ Let $T_A:S_A \to S_A$ be the composition $D_1 \circ D_2 \dots D_n$,
which we call the Coxeter diffeomorphism of the planar tree $A$.

\midinsert
\cline{\epsffile{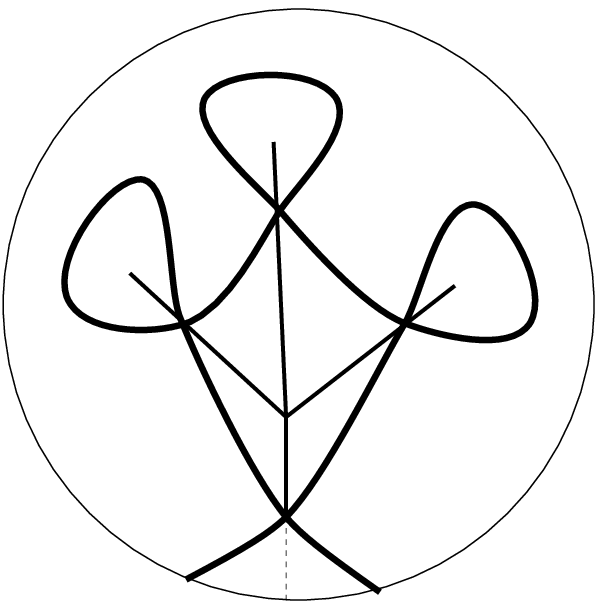}}
\centerline { Figure $3$:
The slalom $Lys$.}
\medskip                 
\endinsert
\goodbreak

The numerical function $A \mapsto \lambda_{H_1}(T_A)$ on trees is
for the inclusion of trees monoton [AC2]. Question: is the function
$A \mapsto \lambda_{\pi_1}(T_A)$  monoton?

\noindent \S 4.{ \bf Trees and hyperbolic knots.}\par

\noindent We now  give the proof of the if part of Theorem $1$:

\demo{\bf Proof}  
Let $B \subset D$ be a rooted tree, such that the Dynkin diagram $\Delta_B$ 
is neither the diagram $A_{2k},\ 1 \leq k,$ nor $E_6$ or $E_8$. We will 
show that
the isotopy class of the monodromy of the fibered knot $K_B$ is
pseudo-anosov.
The geometric 
monodromy $T$ of the knot $K_B$ is up to conjugacy the diffeomorphism  
$T_{A}:S_{A} \to S_{A},$ where we put $A:=\Delta_B$, see [AC4]. The action 
of $T_A$ on
the first homology of $S_A$ is conjugated to the skew Coxeter element $sC_A$
of the tree $A$.
It follows from [AC2] 
that the biggest absolute value $s$ of an eigenvalue of the action of $T_{A}$ 
on the first homology of $S_{A}$ strictly exceeds $1.$ 
So for the homological entropy 
we have $\lambda_{H_1}(T) = \log(s) > 0.$ 
By the entropy inequality, we deduce for the isotopical entropy 
$\lambda_{isotop}(T)$ the inequalities:

$$
0< \lambda_{H_1}(T) \leq \lambda_{\pi_1}(T) \leq \lambda_{isotop}(T) 
\leq \lambda_{top}(T)
$$

\noindent where $\lambda_{isotop}(T)$ is the minimum of the topological entropy 
$\lambda_{top}(T)$ over the relative 
isotopy class of $T.$ Since the isotopical entropy of $T$ is positive, we
conclude that in the decomposition of Thurston [T1]
in quasi-finite and pseudo-anosov pieces of
the diffeomorphism $T$ at least one  pseudo-anosov piece occurs. 
So, to prove that the
isotopy class of the
diffeomorphism $T$ is pseudo-anosov, we need to
prove
that $T$ is irreducible. A reduction of the diffeomorphism $T$ would  give
an essential torus in the complement of the knot $K_B$. Since the knot $K_B$
is an arboricent knot, as shown by the construction of Mikami Hirasawa, we 
conclude with the proposition $2.1$ of [B-Z], see also [O], that the complement of the knot 
$K_B$ does not have an essential torus. So, the diffeomorphism $T$ is
irreducible and hence  pseudo-anosov. 
We can conclude with a celebrated Theorem of W. Thurston [T2], see [O],
that the mapping torus of the diffeomorphism $T$, which is diffeomorphic to
the complement of the knot $K_B$, admits a complete hyperbolic metric.\enddemo

The knot of the slalom curve $E_8$ of the rooted tree with Cayley code 
$[0,1,1,2]$ is not hyperbolic (see Fig. $2$).  The knots of the 
slalom curve Lys of the tree 
with code $[0,1,1,1]$ and of the 
slalom curve $E_{10}$ of the tree $[0,1,1,2,4]$ are hyperbolic 
(see Fig. $1,3$). 

\midinsert
\cline{\epsffile{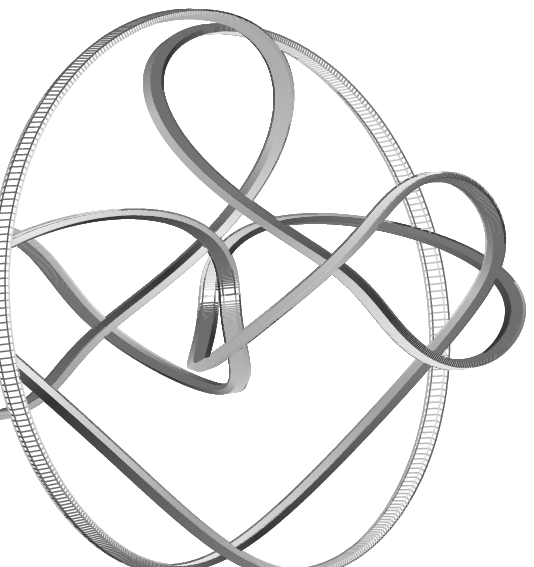}  \epsffile{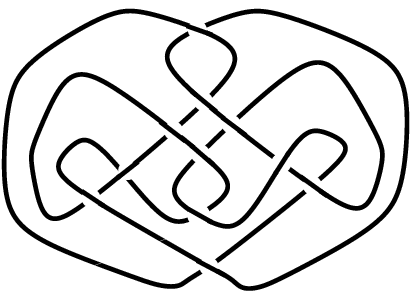}}
\medskip
\centerline{Figure $4$: The knot of the slalom  Lys.}
\endinsert

In fact, for a diffeomorphism $T$ of surfaces the equality
$\lambda_{\pi_1}(T)=\lambda_{isotop}(T)$
holds, and moreover, for pseudo-anosov diffeomorphisms the equality
$\lambda_{\pi_1}(T)=\lambda_{top}(T)$
holds. It would be very interesting to
compute the hyperbolic volume of the knot of the rooted tree $E_{10}$
and to relate it with $\lambda_{isotop}(T_{E_{10}}).$ 

The knots $K_B$ for $B$ such that the Dynkin $\Delta_B$ diagram equals 
$A_{2n},E_6$ or $E_8$, are links of singularities
and the corresponding monodromies are irreducible and of finite order. 
From this fact and from the
proof of the theorem we deduce that Coxeter diffeomorphism of rooted trees 
have in general a pseudo-anosov isotopy class. More precisely, with 
the notations of
section 2 we have:

\proclaim{Theorem $3$}
Let $A \subset D$ be a rooted, bicolored, tree embedded in the plane such that
the root is a terminal vertex and that no other terminal vertex has the color
of the root. The Coxeter
diffeomorphism
$T_A:S_A \to S_A$ is irreducible. Moreover, if $A$ is not
among $A_{2n},E_6,E_8$,
the diffeomorphism is pseudo-anosov.
\endproclaim

The pseudo-anosov diffeomorphisms given by this theorem  
are products of Dehn twists, which belong to
the same conjugacy class in the mapping class group of a surface with one
boundary component and 
the union of
the cores of the Dehn twists of the product decomposition is a spline in the
surface. The pseudo-anosov diffeomorphisms, which we obtain here, differ from
the examples of R. C. Penner [P], see also [F], since all Dehn twist in the
product  belong to the same conjugacy class. 
The diffeomorphism is pseudo-anosov, if
and only if the Dynkin diagram of the intersection of the core curves is not
a classical Dynkin diagram of a finite Coxeter group. A finite 
tree can be realized as Dynkin diagram of a slalom divide of a  (disk wide) web,
which  we define as an embedded finite 
tree $B$ in the
unit disk $D$ such that the intersection $B \cap \partial{D}$ is a set of
terminal vertices of $B$, which are called root vertices of $B$.
The
definition of slalom curve remains unchanged, except for
the slalom of a web
without root vertices, where we consider an immersion of the circle instead
of
the interval. For instance the extended Dynkin diagram $\tilde{E}_8$ with $9$
vertices is the Dynkin diagram of the slalom of Fig. $5$, which is the slalom
of a web with $2$ root
vertices.

\midinsert
\cline{\epsffile{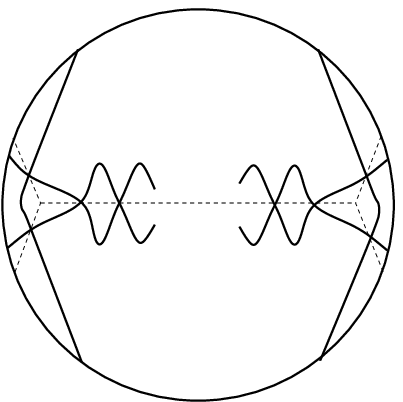}  \epsffile{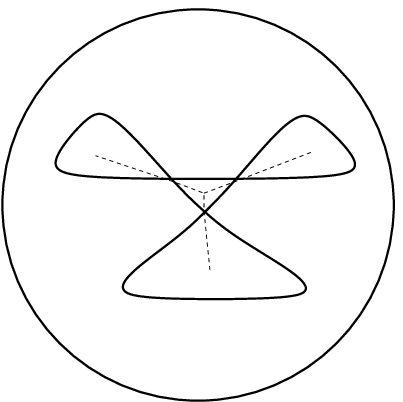}}
\medskip
\cline{\epsffile{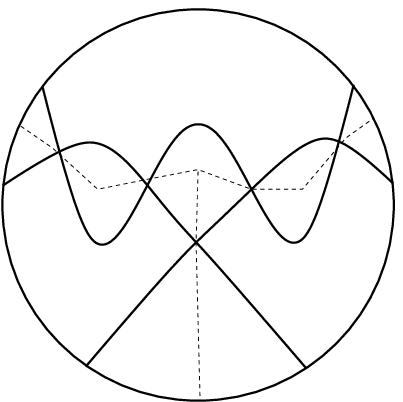}  \epsffile{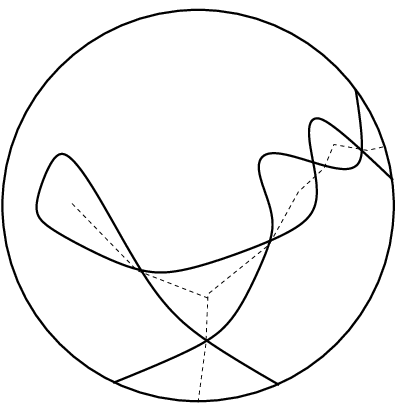}}
\centerline{Figure $5$: The slalom of $\tilde{D}_n,\, n \geq 4$,
$\tilde{E}_6$,$\tilde{E}_7$ and $\tilde{E}_8$.}
\endinsert

The link of
the slalom divide $\tilde{E}_8$ has $2$ components, which is the 
Montesinos link
$M(0,(2,1),(3,1),(6,1))$. The extended Dynkin diagram
$\tilde{D}_4$ with $5$ vertices corresponds to the slalom of the web  with
$4$ root vertices and a single vertex of valency $4$ and its   link  
is the  Montesinos link $M(0,(2,1),(2,1),(2,1),(2,1))$. It
is interesting to observe that both links are in the list $b).$ of proposition
$2.1$ of [B-Z]. The Coxeter diffeomorphism of a Dynkin diagram, which we
suppose to be a tree here, 
is always the monodromy diffeomorphism of a fibered link
by the fibration theorem of [AC4]. The Dynkin diagram $\tilde{D}_4$ is
realized as the
Dynkin diagram of a complete intersection curve singularity by Marc Giusti [G]
and
the Coxeter diffeomorphism of $\tilde{D}_4$ appears as monodromy in the
unfolding of this singularity. The fibered 
link of the web with $n+1$ vertices
and $n$ root vertices, $n \geq 5$, is a chain with $n$ links. This link is
studied in the lecture notes of W. Thurston and is hyperbolic.

\remark {\bf Problem} The complexity $C_{\pi_1}(\phi)$ of an orientation 
preserving isotopy class $\phi$ of 
diffeomorphisms of a 
surface can be defined as the minimum of the quantity $a+b$  over all the
product decompositions of $\phi$ as product of Dehn twists, where $a$ is the
length of the product and where $b$ is the number of intersection points of
the core curves of the  twists involved in the product decomposition. 
The corresponding homological
complexity is the complexity $C_{H_1}(\phi)$ 
where we minimize the quantity $a+h$,
where $h$ stands for the sum of the absolute values of the mutual homological
intersection numbers of the core curves.

The homological complexity of the monodromy $T$ of a non trivial fibered knot 
is estimated from below by $C_{H_1}(T) \geq 4\delta-1$, where 
$\delta$ is the genus of the
fiber. So, we can observe  that 
both complexities coincide and are minimal with respect to this estimation by
the genus  for monodromies of
knots of slalom curves. It would be nice
to  deduce from this observation 
that the homological and isotopical entropy of monodromies of knots of slalom
curves coincide. 
\endremark

\par

\bye